\documentclass[1p,12pt]{elsarticle}

\usepackage{amssymb}
\usepackage{amsmath}
\usepackage{graphicx}

\journal{Applied Mathematics and Computation}
\begin{document}

\begin{frontmatter}

\title{The complete set of homogeneous Hilbert curves in two dimensions}

\author[eer]{C. P\'erez-Demydenko}

\author[baf]{I. Brito-Reyes}

\author[baf]{B. Arag\'on Fern\'andez}

\author[eer]{E. Estevez-Rams\corref{cor1}}
\ead{estevez@imre.oc.uh.cu}

\cortext[cor1]{Corresponding author at: Instituto de Ciencias y Tecnolog\'ia de Materiales, University of Havana (IMRE), San Lazaro y L. CP 10400. La Habana. Cuba. t: (+537) 8705707}

\address[eer]{Instituto de Ciencias y Tecnolog\'ia de Materiales, University of Havana (IMRE), San Lazaro y L. CP 10400. La Habana. Cuba.}

\address[baf]{Universidad de las Ciencias Inform\'aticas (UCI), Carretera a San Antonio, Boyeros. La Habana. Cuba.}

\begin{abstract}
An inventory of all possible homogenous Hilbert curves in two dimensions are reported. Six new Hilbert curves are described by introducing the reversion operation in the construction algorithm. For each curve, the set of affine transformation defining the generation process is reported. Finally, each curve is also described in terms of a tag system.
\end{abstract}

\begin{keyword}
Hilbert curve \sep Space filling curves  
\end{keyword}


\date{\today}
\end{frontmatter}

\section{Introduction}

A space filling curve is a continuous curve that covers completely the d-dimensional space. In mathematical terms, is a surjective mapping from  $\mathbb{R}\longrightarrow \mathbb{R}^{d}$. The concept was described by Peano and shortly after, Hilbert introduced, in two dimensions (2D), an optimized space filling curve that now bears his name \cite{sagan94}.

Hilbert curves have found applications in several fields including, computer science \cite{chen05,chen11,bially69}, manufacturing engineering \cite{cox94}, imaging processing \cite{songa02,liang08}, data visualization \cite{anders09}, and the traveling salesman problem \cite{gao94}, among others. The wealth of applications stems from the good locality preserving properties of such curves, meaning that neighbors interval in the one dimensional curve, are fairly well clustered in the $\mathbb{R}^{d}$ mapping \cite{moon01,bauman06}.

In spite that the original Hilbert construction is dealt almost exclusively, there are several Hilbert curves for a given dimension. Haverkort \cite{haverkort11} discussed the set of Hilbert curves in three dimensions, finding more than $10^{6}$ possible variants of Hilbert curves. In 2D, already Moore described a curve that is named after him \cite{moore00}. More recently, Liu  has discussed four alternative patterns of the Hilbert curve in 2D \cite{liu04}. Yet, Liu description is somewhat confusing, and the description of the generating process  incomplete.

In this contribution the complete set of homogeneous Hilbert curves in 2D are reported. In doing so, six new patterns are introduced. The recursive, analytical, generation of all twelve patterns is discussed, which allows to follow the logic of the construction for such curves. The Hilbert curves are also described in terms of a tag-system.

\section{Hilbert construction}\label{sec:hilbconst}

For a given $n$, the unit segment $[0,1]$ is divided into $4^{n}$ non-overlapping intervals $\{I_{i}\}_{n}$ of equal length, that completely covers the original segment. Correspondingly, the unit square $[0,1] \times [0,1]$ is partitioned into $4^{n}$ non-overlapping squares $\{S_{i}\}_{n}$ of equal area, that completely covers the original square. A mapping can be set from $\{I_{i}\}_{n}$ onto $\{S_{i}\}_{n}$.

We proceed as follows: For $n+1$ the unit segment $[0,1]$ is divided in such a way, that each $I_{j} \epsilon \{I_{i}\}_{n}$ is partitioned into four equal length segments $\{I'_{j}, I'_{j+1}, I'_{j+2}, I'_{j+3}\}$, each of them belonging to $\{I_{i}\}_{n+1}$. Similarly, each $S_{k} \epsilon \{S_{i}\}_{n}$ is partitioned into four equal squares  $\{S'_{j}, S'_{j+1}, S'_{j+2}, S'_{j+3}\}$, each of them  belonging to $\{S_{i}\}_{n+1}$. In  this way results an infinite sequence of $\{I_{i}\}_{n}$ $(n=1,2, \dots)$ partitions of the unit segment, and a corresponding sequence of $\{S_{j}\}_{n}$ $(n=1,2, \dots)$ partitions of the unit square. The following two conditions are imposed to the mapping from $\{I_{j}\}_{n}$ onto $\{S_{j}\}_{n}$ (Hilbert conditions):

\begin{enumerate}
 \item For a fixed $n$, two adjacent intervals $I_{j}, I_{j+1} \epsilon \{I_{i}\}_{n}$ will be mapped to two adjacent squares (sharing an edge) $S_{k}, S_{k+1} \epsilon \{S_{i}\}_{n}$.

\item Taken $I_{j} \epsilon \{I_{i}\}_{n}$, let $S_{k} \epsilon \{S_{i}\}_{n}$ be the corresponding square under the mapping then, if $I'_{l} \epsilon \{I_{i}\}_{n+1}$ is one of the intervals in the partition of $I_{j}$, its image $S'_{p} \epsilon \{S_{i}\}_{n+1}$ under the mapping, will be one of the squares in the partition of $S_{k}$.
\end{enumerate}

Hilbert conditions do not uniquely determine the mapping between $\{I_{i}\}_{n}$ and $\{S_{i}\}_{n}$ \cite{moore00}. In order to fix uniquely the mapping, boundary conditions must be imposed, by stating to which squares $S_{j}$ and $S_{k}$, will the initial $I_{0}$ and final $I_{4^{n}-1}$ interval be mapped to, respectively.

Hilbert curve will then be the space filling curve resulting from the limit process $n \longrightarrow \infty$. The mapping $\{I_{i}\}_{n}\longrightarrow \{S_{i}\}_{n}$ is also known as the Hilbert curve approximation of order $n$, the term Hilbert curve will be used for this approximation of order $n$, as the context will make its use unambiguous. 

In what follows, all Hilbert curves related by rotation, mirror or reversal operations will be considered equivalent.

\section{The six proper Hilbert curves}\label{sec:sixhilbcurves}

As the unit square is divided into quadrants, it is needed to specify the connectivity between quadrants. Quadrants will be numbered in clockwise order, starting from the left lower quadrant, and connectivity will follow the same order: $0 \longrightarrow 1 \longrightarrow 2 \longrightarrow 3$. There is no connectivity between quadrant $0$ and $3$ (Figure \ref{fig:quadrants}).

For $n=1$, there is only one possible Hilbert curve $H_{1}$, depicted in figure \ref{fig:h1}. The circle flags the starting or entry point of the curve, and the arrow head, the ending or exit point of the curve.

For $n=2$, Hilbert's conditions result in the curves shown in figure \ref{fig:h2}. Curve \ref{fig:h2}a1 is the original Hilbert curve, which will be denoted by $_{0}H_{n}$; curve \ref{fig:h2}b1 is the one described by Moore, which will be denoted by $_{1}H_{n}$; while curves \ref{fig:h2}a2 to \ref{fig:h2}b4, are the four curves described by Liu \cite{liu04}, and will be numbered from $_{3}H_{n}$ to $_{5}H_{n}$ in row scan order. 

It can be seen, that the six curves in figure \ref{fig:h2} have in each quadrant, a $\frac{1}{2}$-scaled $H_{1}$ curve, rotated or mirror reflected, If the $H_{n+1}$ curve is described by the rotation or mirror operation, carried out to the scaled $H_{n}$ for each quadrant then, the six curves $_{i}H_{2}$ can be depicted as shown on figure \ref{fig:h2op}. It follows, that constructing $H_{n+1}$ from $H_{n}$ involves an affine transformation $q_{i}$ of $H_{n}$ for each quadrant $i$. The affine transformation $q_{i}$ has a scaling rotational part described by a $2 \times 2$ matrix, and a translation given by a vector.

For the original Hilbert curve $_{0}H_{n}$ (figure \ref{fig:h2}a1), the affine transformation will be given by

\begin{equation*}
\begin{array}{llll}
_{0}q_{0}\left ( \begin{array}{c} x \\ y \end{array}\right )= & \frac{1}{2}\left ( \begin{array}{cc} 0 & 1  \\ 1 & 0 \end{array} \right )\left ( \begin{array}{c} x \\ y \end{array}\right ) & + & \frac{1}{2}\left ( \begin{array}{c} 0 \\ 0 \end{array}\right ) \\
_{0}q_{1}\left ( \begin{array}{c} x \\ y \end{array}\right )= & \frac{1}{2}\left ( \begin{array}{cc} 1 & 0  \\ 0 & 1 \end{array} \right )\left ( \begin{array}{c} x \\ y \end{array}\right ) & + & \frac{1}{2}\left ( \begin{array}{c} 0 \\ 1 \end{array}\right ) \\ 
_{0}q_{2}\left ( \begin{array}{c} x \\ y \end{array}\right )= & \frac{1}{2}\left ( \begin{array}{cc} 1 & 0  \\ 0 & 1 \end{array} \right )\left ( \begin{array}{c} x \\ y \end{array}\right ) & + & \frac{1}{2}\left ( \begin{array}{c} 1 \\ 1 \end{array}\right ) \\
_{0}q_{3}\left ( \begin{array}{c} x \\ y \end{array}\right )= & \frac{1}{2}\left ( \begin{array}{cc} 0 & -1 \\ -1 & 0\end{array} \right )\left ( \begin{array}{c} x \\ y \end{array}\right ) & + & \frac{1}{2}\left ( \begin{array}{c} 2 \\ 1 \end{array}\right ) .
\end{array}
\end{equation*}

Where $q_{i}$ transform the $H_{n}$ curve into the curve in the $i$-quadrant of the $H_{n+1}$ mapping.

For the Moore curve $_{1}H_{n}$ (figure \ref{fig:h2}b1), the affine transformation will be given by

\begin{equation*}
\begin{array}{llll}
_{1}q_{0}\left ( \begin{array}{c} x \\ y \end{array}\right )= & \frac{1}{2}\left ( \begin{array}{cc} 0 & -1 \\ 1 & 0 \end{array} \right )\left ( \begin{array}{c} x \\ y \end{array}\right ) & + & \frac{1}{2}\left ( \begin{array}{c}  1 \\ 0 \end{array}\right ) \\
_{1}q_{1}\left ( \begin{array}{c} x \\ y \end{array}\right )= & \frac{1}{2}\left ( \begin{array}{cc} 0 & -1 \\  1 & 0 \end{array} \right )\left ( \begin{array}{c} x \\ y \end{array}\right ) & + & \frac{1}{2}\left ( \begin{array}{c} 1 \\ 1 \end{array}\right ) \\ 
_{1}q_{2}\left ( \begin{array}{c} x \\ y \end{array}\right )= & \frac{1}{2}\left ( \begin{array}{cc} 0 & 1  \\ -1 & 0 \end{array} \right )\left ( \begin{array}{c} x \\ y \end{array}\right ) & + & \frac{1}{2}\left ( \begin{array}{c} 1 \\ 2 \end{array}\right ) \\
_{1}q_{3}\left ( \begin{array}{c} x \\ y \end{array}\right )= & \frac{1}{2}\left ( \begin{array}{cc} 0 & 1  \\ -1 & 0 \end{array} \right )\left ( \begin{array}{c} x \\ y \end{array}\right ) & + & \frac{1}{2}\left ( \begin{array}{c} 1 \\ 1 \end{array}\right ). \\
\end{array}
\end{equation*}

The affine transformation for the four curves described by Liu \cite{liu04}, will be

\begin{equation*}
\begin{array}{llll}
_{2}q_{0}\left ( \begin{array}{c} x \\ y \end{array}\right )= & \frac{1}{2}\left ( \begin{array}{cc} -1 & 0 \\  0 & -1 \end{array} \right )\left ( \begin{array}{c} x \\ y \end{array}\right ) & + & \frac{1}{2}\left ( \begin{array}{c} 1 \\ 1 \end{array}\right ) \\
_{2}q_{1}\left ( \begin{array}{c} x \\ y \end{array}\right )= & \frac{1}{2}\left ( \begin{array}{cc}  1 & 0 \\  0 & 1 \end{array} \right )\left ( \begin{array}{c} x \\ y \end{array}\right ) & + & \frac{1}{2}\left ( \begin{array}{c}  0 \\ 1 \end{array}\right ) \\ 
_{2}q_{2}\left ( \begin{array}{c} x \\ y \end{array}\right )= & \frac{1}{2}\left ( \begin{array}{cc}  1 & 0 \\  0 & 1 \end{array} \right )\left ( \begin{array}{c} x \\ y \end{array}\right ) & + & \frac{1}{2}\left ( \begin{array}{c}  1 \\ 1 \end{array}\right ) \\
_{2}q_{3}\left ( \begin{array}{c} x \\ y \end{array}\right )= & \frac{1}{2}\left ( \begin{array}{cc} -1 & 0 \\  0 & -1 \end{array} \right )\left ( \begin{array}{c} x \\ y \end{array}\right ) & + & \frac{1}{2}\left ( \begin{array}{c} 2 \\ 1 \end{array}\right ). \\
\end{array}
\end{equation*}
 
for $_{2}H_{n}$ (figure \ref{fig:h2}a2);

\begin{equation*}
\begin{array}{llll}
_{3}q_{0}\left ( \begin{array}{c} x \\ y \end{array}\right )= & \frac{1}{2}\left ( \begin{array}{cc} 1 & 0  \\  0 & -1 \end{array} \right )\left ( \begin{array}{c} x \\ y \end{array}\right ) & + & \frac{1}{2}\left ( \begin{array}{c} 0 \\ 1 \end{array}\right ) \\
_{3}q_{1}\left ( \begin{array}{c} x \\ y \end{array}\right )= & \frac{1}{2}\left ( \begin{array}{cc} 0 & -1 \\  1 & 0 \end{array} \right )\left ( \begin{array}{c} x \\ y \end{array}\right ) & + & \frac{1}{2}\left ( \begin{array}{c} 1 \\ 1 \end{array}\right ) \\ 
_{3}q_{2}\left ( \begin{array}{c} x \\ y \end{array}\right )= & \frac{1}{2}\left ( \begin{array}{cc} 0 & 1  \\  -1 & 0 \end{array} \right )\left ( \begin{array}{c} x \\ y \end{array}\right ) & + & \frac{1}{2}\left ( \begin{array}{c} 1 \\ 2 \end{array}\right ) \\
_{3}q_{3}\left ( \begin{array}{c} x \\ y \end{array}\right )= & \frac{1}{2}\left ( \begin{array}{cc} 1 & 0  \\  0 & -1 \end{array} \right )\left ( \begin{array}{c} x \\ y \end{array}\right ) & + & \frac{1}{2}\left ( \begin{array}{c} 1 \\ 1 \end{array}\right ),\\
\end{array}
\end{equation*}

for $_{3}H_{n}$ (figure \ref{fig:h2}b2);

\begin{equation*}
\begin{array}{llll}
_{4}q_{0}\left ( \begin{array}{c} x \\ y \end{array}\right )= & \frac{1}{2}\left ( \begin{array}{cc} 0 & 1  \\  1 & 0 \end{array} \right )\left ( \begin{array}{c} x \\ y \end{array}\right ) & + & \frac{1}{2}\left ( \begin{array}{c} 0 \\ 0 \end{array}\right ) \\
_{4}q_{1}\left ( \begin{array}{c} x \\ y \end{array}\right )= & \frac{1}{2}\left ( \begin{array}{cc} 1 & 0  \\  0 & 1 \end{array} \right )\left ( \begin{array}{c} x \\ y \end{array}\right ) & + & \frac{1}{2}\left ( \begin{array}{c} 0 \\ 1 \end{array}\right ) \\ 
_{4}q_{2}\left ( \begin{array}{c} x \\ y \end{array}\right )= & \frac{1}{2}\left ( \begin{array}{cc} 1 & 0  \\  0 & 1 \end{array} \right )\left ( \begin{array}{c} x \\ y \end{array}\right ) & + & \frac{1}{2}\left ( \begin{array}{c} 1 \\ 1 \end{array}\right ) \\
_{4}q_{3}\left ( \begin{array}{c} x \\ y \end{array}\right )= & \frac{1}{2}\left ( \begin{array}{cc} -1 & 0 \\  0 & -1 \end{array} \right )\left ( \begin{array}{c} x \\ y \end{array}\right ) & + & \frac{1}{2}\left ( \begin{array}{c} 2 \\ 1 \end{array}\right ), \\
\end{array}
\end{equation*}

for $_{4}H_{n}$ (figure \ref{fig:h2}a3);

\begin{equation*}
\begin{array}{llll}
_{5}q_{0}\left ( \begin{array}{c} x \\ y \end{array}\right )= & \frac{1}{2}\left ( \begin{array}{cc} 1 & 0  \\  0 & -1 \end{array} \right )\left ( \begin{array}{c} x \\ y \end{array}\right ) & + & \frac{1}{2}\left ( \begin{array}{c} 0 \\ 1 \end{array}\right ) \\
_{5}q_{1}\left ( \begin{array}{c} x \\ y \end{array}\right )= & \frac{1}{2}\left ( \begin{array}{cc} 0 & -1 \\  1 & 0  \end{array} \right )\left ( \begin{array}{c} x \\ y \end{array}\right ) & + & \frac{1}{2}\left ( \begin{array}{c} 1 \\ 1 \end{array}\right ) \\ 
_{5}q_{2}\left ( \begin{array}{c} x \\ y \end{array}\right )= & \frac{1}{2}\left ( \begin{array}{cc} 0 & 1  \\  -1 & 0 \end{array} \right )\left ( \begin{array}{c} x \\ y \end{array}\right ) & + & \frac{1}{2}\left ( \begin{array}{c} 1 \\ 2 \end{array}\right ) \\
_{5}q_{3}\left ( \begin{array}{c} x \\ y \end{array}\right )= & \frac{1}{2}\left ( \begin{array}{cc} 0 & 1  \\  -1 & 0 \end{array} \right )\left ( \begin{array}{c} x \\ y \end{array}\right ) & + & \frac{1}{2}\left ( \begin{array}{c} 1 \\ 1 \end{array}\right ). \\
\end{array}
\end{equation*}

for $_{5}H_{n}$ (figure \ref{fig:h2}b3).

It could be considered straight forward to apply recursively the corresponding affine transformation of each Hilbert curve to $_{k}H_{n}$, in order to obtain $_{k}H_{n+1}$. The problem is, that except for $_{0}H_{n}$, applying $_{k}q$ to $_{k}H_{2}$ will render impossible to achieve the desired quadrant connectivity for $k=1,2,3,4,5$ (a point missed in Liu description \cite{liu04}). Instead, to build $_{k}H_{n}$, $_{0}H_{n-1}$ is constructed by applying recursively $_{0}q$ and then,  $_{k}H_{n}$ is built, by applying $_{k}q$ to $_{0}H_{n-1}$. The result for $n=3$ and $n=4$ are shown in figure \ref{fig:h3} and figure \ref{fig:h4}.

\section{The six improper homogeneous Hilbert curves}\label{sec:hc}

By construction, the only Hilbert curves $_{k}H_{n}$ that can be used to generate a higher order curve $_{k}H_{n+1}$, are those with initial and ending point at adjacent corners or edges. Inspection of figure \ref{fig:h2} shows, that only $_{0}H_{n}$ and $_{5}H_{n}$ are candidate generating curves. $_{0}H_{n}$ was already used for the construction of the first six curves. 

As already said, for all curves, other than $_{0}H_{n}$, the recursive application of their corresponding affine mapping, does not yield a well connected $_{k}H_{n+1}$ curve. In order to construct a higher order Hilbert curve from $_{5}H_{n}$, a new operation called  reversion is introduced, which exchange the entry and exit point of the $_{k}H_{n}$ curve, as described in figure \ref{fig:swap}. 

Six new Hilbert curves can be constructed using the reversion operation. Such curves will be called improper Hilbert curves, to differentiate them from the earlier six curves. In terms of the affine transformation performed in each quadrant, the new  curves are described in figure \ref{fig:handimproper}. The curves for $n=3$ and $n=4$ are shown in figures \ref{fig:improper3} and \ref{fig:improper4}, respectively. The affine transformation for the first curve $_{6}H_{n}$ will be

\begin{equation*}
\begin{array}{llll}
_{6}q_{0}\left ( \begin{array}{c} x \\ y \end{array}\right )= & \frac{1}{2}\left ( \begin{array}{cc} -1 & 0  \\  0 & -1 \end{array} \right )\left ( \begin{array}{c} x \\ y \end{array}\right ) & + & \frac{1}{2}\left ( \begin{array}{c} 1 \\ 1 \end{array}\right ) \\
_{6}\bar{q}_{1}\left ( \begin{array}{c} x \\ y \end{array}\right )= & \frac{1}{2}\left ( \begin{array}{cc} -1 & 0 \\  0 & 1  \end{array} \right )\left ( \begin{array}{c} x \\ y \end{array}\right ) & + & \frac{1}{2}\left ( \begin{array}{c} 1 \\ 1 \end{array}\right ) \\ 
_{6}q_{2}\left ( \begin{array}{c} x \\ y \end{array}\right )= & \frac{1}{2}\left ( \begin{array}{cc} 1 & 0  \\  0 & 1 \end{array} \right )\left ( \begin{array}{c} x \\ y \end{array}\right ) & + & \frac{1}{2}\left ( \begin{array}{c} 1 \\ 1 \end{array}\right ) \\
_{6}\bar{q}_{3}\left ( \begin{array}{c} x \\ y \end{array}\right )= & \frac{1}{2}\left ( \begin{array}{cc} 1 & 0  \\  0 & -1 \end{array} \right )\left ( \begin{array}{c} x \\ y \end{array}\right ) & + & \frac{1}{2}\left ( \begin{array}{c} 1 \\ 1 \end{array}\right ). \\
\end{array}
\end{equation*}

but now the bar over the $_{6}\bar{q}_{1}$ and $_{6}\bar{q}_{3}$ means, that besides the affine transformation, the reversion operation also acts over $_{5}H_{n}$. For the second improper Hilbert curve $_{7}H_{n}$, the affine transformation will be

\begin{equation*}
\begin{array}{llll}
_{7}q_{0}\left ( \begin{array}{c} x \\ y \end{array}\right )= & \frac{1}{2}\left ( \begin{array}{cc} -1 & 0  \\  0 & -1 \end{array} \right )\left ( \begin{array}{c} x \\ y \end{array}\right ) & + & \frac{1}{2}\left ( \begin{array}{c} 1 \\ 1 \end{array}\right ) \\
_{7}\bar{q}_{1}\left ( \begin{array}{c} x \\ y \end{array}\right )= & \frac{1}{2}\left ( \begin{array}{cc} -1 & 0 \\  0 & 1  \end{array} \right )\left ( \begin{array}{c} x \\ y \end{array}\right ) & + & \frac{1}{2}\left ( \begin{array}{c} 1 \\ 1 \end{array}\right ) \\ 
_{7}q_{2}\left ( \begin{array}{c} x \\ y \end{array}\right )= & \frac{1}{2}\left ( \begin{array}{cc} 1 & 0  \\  0 & 1 \end{array} \right )\left ( \begin{array}{c} x \\ y \end{array}\right ) & + & \frac{1}{2}\left ( \begin{array}{c} 1 \\ 1 \end{array}\right ) \\
_{7}q_{3}\left ( \begin{array}{c} x \\ y \end{array}\right )= & \frac{1}{2}\left ( \begin{array}{cc} 0 & -1 \\ -1 & 0\end{array} \right )\left ( \begin{array}{c} x \\ y \end{array}\right ) & + & \frac{1}{2}\left ( \begin{array}{c} 2 \\ 1 \end{array}\right ). \\
\end{array}
\end{equation*}

For the third improper Hilbert curve $_{8}H_{n}$, the affine transformation will be

\begin{equation*}
\begin{array}{llll}
_{8}\bar{q}_{0}\left ( \begin{array}{c} x \\ y \end{array}\right )= & \frac{1}{2}\left ( \begin{array}{cc} 0 & 1  \\  -1 & 0 \end{array} \right )\left ( \begin{array}{c} x \\ y \end{array}\right ) & + & \frac{1}{2}\left ( \begin{array}{c} 0 \\ 1 \end{array}\right ) \\
_{8}\bar{q}_{1}\left ( \begin{array}{c} x \\ y \end{array}\right )= & \frac{1}{2}\left ( \begin{array}{cc} -1 & 0 \\  0 & 1  \end{array} \right )\left ( \begin{array}{c} x \\ y \end{array}\right ) & + & \frac{1}{2}\left ( \begin{array}{c} 1 \\ 1 \end{array}\right ) \\ 
_{8}q_{2}\left ( \begin{array}{c} x \\ y \end{array}\right )= & \frac{1}{2}\left ( \begin{array}{cc} 1 & 0  \\  0 & 1 \end{array} \right )\left ( \begin{array}{c} x \\ y \end{array}\right ) & + & \frac{1}{2}\left ( \begin{array}{c} 1 \\ 1 \end{array}\right ) \\
_{8}q_{3}\left ( \begin{array}{c} x \\ y \end{array}\right )= & \frac{1}{2}\left ( \begin{array}{cc} 0 & -1 \\ -1 & 0\end{array} \right )\left ( \begin{array}{c} x \\ y \end{array}\right ) & + & \frac{1}{2}\left ( \begin{array}{c} 2 \\ 1 \end{array}\right ). \\
\end{array}
\end{equation*}

For the fourth improper Hilbert curve $_{9}H_{n}$, the affine transformation will be

\begin{equation*}
\begin{array}{llll}
_{9}\bar{q}_{0}\left ( \begin{array}{c} x \\ y \end{array}\right )= & \frac{1}{2}\left ( \begin{array}{cc} 0 & -1  \\  -1 & 0 \end{array} \right )\left ( \begin{array}{c} x \\ y \end{array}\right ) & + & \frac{1}{2}\left ( \begin{array}{c} 1 \\ 1 \end{array}\right ) \\
_{9}q_{1}\left ( \begin{array}{c} x \\ y \end{array}\right )= & \frac{1}{2}\left ( \begin{array}{cc} 0 & -1 \\  1 & 0  \end{array} \right )\left ( \begin{array}{c} x \\ y \end{array}\right ) & + & \frac{1}{2}\left ( \begin{array}{c} 1 \\ 1 \end{array}\right ) \\ 
_{9}\bar{q}_{2}\left ( \begin{array}{c} x \\ y \end{array}\right )= & \frac{1}{2}\left ( \begin{array}{cc} 0 & 1  \\  1 & 0 \end{array} \right )\left ( \begin{array}{c} x \\ y \end{array}\right ) & + & \frac{1}{2}\left ( \begin{array}{c} 1 \\ 1 \end{array}\right ) \\
_{9}q_{3}\left ( \begin{array}{c} x \\ y \end{array}\right )= & \frac{1}{2}\left ( \begin{array}{cc} 0 & 1 \\ -1 & 0\end{array} \right )\left ( \begin{array}{c} x \\ y \end{array}\right ) & + & \frac{1}{2}\left ( \begin{array}{c} 1 \\ 1 \end{array}\right ). \\
\end{array}
\end{equation*}

For the fifth improper Hilbert curve $_{10}H_{n}$, the affine transformation will be

\begin{equation*}
\begin{array}{llll}
_{10}q_{0}\left ( \begin{array}{c} x \\ y \end{array}\right )= & \frac{1}{2}\left ( \begin{array}{cc} 1 & 0  \\  0 & -1 \end{array} \right )\left ( \begin{array}{c} x \\ y \end{array}\right ) & + & \frac{1}{2}\left ( \begin{array}{c} 0 \\ 1 \end{array}\right ) \\
_{10}q_{1}\left ( \begin{array}{c} x \\ y \end{array}\right )= & \frac{1}{2}\left ( \begin{array}{cc} 0 & -1 \\  1 & 0  \end{array} \right )\left ( \begin{array}{c} x \\ y \end{array}\right ) & + & \frac{1}{2}\left ( \begin{array}{c} 1 \\ 1 \end{array}\right ) \\ 
_{10}\bar{q}_{2}\left ( \begin{array}{c} x \\ y \end{array}\right )= & \frac{1}{2}\left ( \begin{array}{cc} 0 & 1  \\  1 & 0 \end{array} \right )\left ( \begin{array}{c} x \\ y \end{array}\right ) & + & \frac{1}{2}\left ( \begin{array}{c} 1 \\ 1 \end{array}\right ) \\
_{10}\bar{q}_{3}\left ( \begin{array}{c} x \\ y \end{array}\right )= & \frac{1}{2}\left ( \begin{array}{cc} -1 & 0 \\ 0 & -1\end{array} \right )\left ( \begin{array}{c} x \\ y \end{array}\right ) & + & \frac{1}{2}\left ( \begin{array}{c} 2 \\ 1 \end{array}\right ). \\
\end{array}
\end{equation*}

For the last improper Hilbert curve $_{11}H_{n}$, the affine transformation will be

\begin{equation*}
\begin{array}{llll}
_{11}q_{0}\left ( \begin{array}{c} x \\ y \end{array}\right )= & \frac{1}{2}\left ( \begin{array}{cc} 1 & 0  \\  0 & -1 \end{array} \right )\left ( \begin{array}{c} x \\ y \end{array}\right ) & + & \frac{1}{2}\left ( \begin{array}{c} 0 \\ 1 \end{array}\right ) \\
_{11}q_{1}\left ( \begin{array}{c} x \\ y \end{array}\right )= & \frac{1}{2}\left ( \begin{array}{cc} 0 & -1 \\  1 & 0  \end{array} \right )\left ( \begin{array}{c} x \\ y \end{array}\right ) & + & \frac{1}{2}\left ( \begin{array}{c} 1 \\ 1 \end{array}\right ) \\ 
_{11}\bar{q}_{2}\left ( \begin{array}{c} x \\ y \end{array}\right )= & \frac{1}{2}\left ( \begin{array}{cc} 0 & 1  \\  1 & 0 \end{array} \right )\left ( \begin{array}{c} x \\ y \end{array}\right ) & + & \frac{1}{2}\left ( \begin{array}{c} 1 \\ 1 \end{array}\right ) \\
_{11}q_{3}\left ( \begin{array}{c} x \\ y \end{array}\right )= & \frac{1}{2}\left ( \begin{array}{cc} 0 & 1 \\ -1 & 0\end{array} \right )\left ( \begin{array}{c} x \\ y \end{array}\right ) & + & \frac{1}{2}\left ( \begin{array}{c} 1 \\ 1 \end{array}\right ). \\
\end{array}
\end{equation*}

Hilbert curves $_{0}H_{n}$, $_{1}H_{n}$, $_{2}H_{n}$, $_{3}H_{n}$, $_{6}H_{n}$, $_{8}H_{n}$, $_{9}H_{n}$ and $_{10}H_{n}$ are mirror symmetric with respect to a vertical line at half the unit square . $_{1}H_{n}$, $_{2}H_{n}$, $_{6}H_{n}$ and $_{9}H_{n}$ are closed curves in the sense of having the entry and exit point in adjacent squares, such curves are useful, for example,  in the traveling salesman problem \cite{bartholdi82}.

\section{Completeness of the homogeneous Hilbert curves}\label{sec:complete}

As already stated, the only building blocks for higher order Hilbert curves are those curves with entry and exit points in different adjacent edges or corners. The reason for that, is the way quadrants are visited (Figure \ref{fig:quadrants}). Continuity of the Hilbert curve, together with the connectivity scheme $1 \longrightarrow 2\longrightarrow 3\longrightarrow 4$ impose constrains on the possible building block and their orientation in each quadrant. As already seen, only Hilbert original curve $_{0}H_{n}$ (Figure \ref{fig:h2}-1a) and $_{5}H_{n}$ (Figure \ref{fig:h2}-3b) can be used for that purpose. Hilbert curves are uniquely determined by their boundary conditions, and each curve preserves their boundary condition for any order \cite{sagan94,moore00}. Sagan has explained the approximating polygonal construction to depict Hilbert curves \cite{sagan94}. For our purpose, it will suffice to say that in each quadrant, the entry and exit point can be represented as a vector going from the former the latter. Such vector will be called the boundary vector (bvector). In such a way, uniqueness of the space filling curve can be described by just showing the bvector at each quadrant, regardless of the curve order. As already stated, two Hilbert curves are considered equivalent if they are related by a rotation, reflection or reversion operation, correspondingly,  bvectors diagrams are equivalent if they are related by any of the described operations. Therefore, a set of equivalent bvector diagrams gives rise to only one distinctive Hilbert curve.

Figure \ref{fig:boundaryvector}, shows the bvectors for $_{0}H$ and $_5{H}$. Homogeneous curves have been defined as those built from just one building block. Homogeneous Hilbert curve of any order, can only be described by suitable rotations, reflections or both, of the corresponding boundary vector at each quadrant while preserving continuity between connected quadrants. 

For the proper curves there are only four possible orientation for the boundary vector in the first quadrant which results in a curve connected to the second quadrant (Figure \ref{fig:propercomplete}-I). Once the boundary vector for the first quadrant is fixed, the boundary vectors for the second and third quadrants are also fixed by the quadrant connectivity (Figure \ref{fig:propercomplete}-II, III). For the fourth quadrant there are two choices for each case, as there is no connectivity constrain at the ending point and this results in the eight boundary vectors diagram of Figure \ref{fig:propercomplete}-IV. The diagrams 1a and 3b form an ''enantiomorphic'' pair that give rise to two equivalent Hilbert curves related by a mirror-reversion operation, the same goes for the ''enantiomorphic'' pair 2b-4b. Both enantiomorphic pairs contribute with one non-equivalent Hilbert curve. 

For the improper curves the analysis follow the same lines. Again there are four possibles boundary vectors for the first quadrant (Figure \ref{fig:impropercomplete}-I), which fix the boundary vector for the second and third quadrant (Figure \ref{fig:impropercomplete}-II, III).  The liberty of choice in the fourth quadrant leads to a total of eight improper boundary vector diagrams (Figure \ref{fig:impropercomplete}-IV). Diagrams pairs 1b-2b and 3b-4b are mirror symmetric and therefore lead to equivalent ''enantiomorphic''  Hilbert curves. 

There are thus, a total of six non-equivalent improper Hilbert curves.

All possible non equivalent boundary vector diagrams have been exhausted and therefore, all possible constructions for the proper and improper Hilbert curves.  It is in this sense that it is claimed that the described twelve curves (six proper and six improper), exhaust all possible homogeneous Hilbert curves in two dimensions.

It must be pointed out that new Hilbert curves can be built if the homogeneity condition is dropped, and curves of order $n$ are allowed to mix $_{0}H_{n-1}$ and $_{5}H_{n-1}$ curves at different quadrants. In fact, in such case, up to $40$ different Hilbert curves in 2D can be built, results will be published elsewhere.
 
\section{Hilbert curves  as a tag system}\label{sec:tag}

Hilbert curves can be viewed as drawings, where the point is never removed from the surface of the paper, made of succession of unit length lines that can go up, down, right or left. This allows to define an alphabet $\Sigma=\{u, d, r, l\}$, where $u$ stands for \textit{up}, $d$ for \textit{down}, $r$ for \textit{right} and $l$ for \textit{left}. A given Hilbert curve $H_{n}$ can now be encoded by a finite length word, with symbols taken from $\Sigma$, describing the pen stroke at each step. For example, $_{0}H_{1}$ will be encoded by the word $urd$, while $_{0}H_{2}$ will be given by $ruluurdrurddldr$.

The tag system describing the Hilbert construction is then, the rules that transforms a word $h_{n}$, encoding $H_{n}$, to the word $h_{n+1}$ that describes $H_{n+1}$. For the original Hilbert curve $_{0}H_{n}$, Seebold \cite{seebold07} reported the following tag system:

\begin{equation}
_{0}h_{n+1}=\delta_{o} (_{0}h_{n}) \; u \; _{0}h_{n} \; r \; _{0}h_{n} \; d \; \delta_{a} (_{0}h_{n})    
\end{equation}

where $\delta_{o}$ and $\delta_{a}$ are the morphism defined as

\begin{equation}
 \begin{array}{llll}
  \delta_{o}(u)=r & \delta_{o}(r)=u & \delta_{o}(d)=l & \delta_{o}(l)=d \\
  \delta_{a}(u)=l & \delta_{a}(r)=d & \delta_{a}(d)=r & \delta_{a}(l)=u.
 \end{array}
\end{equation}

Similarly, a tag system can be defined for the rest of the proper Hilbert curves as

\begin{equation}
\begin{array}{l}
_{1}h_{n+1}=\delta_{g}(_{0}h_{n})  \; u \; \delta_{g}(_{0}h_{n})  \; r \; \delta_{x}(_{0}h_{n}) \; d \; \delta_{x}(_{0}h_{n}) \\
_{2}h_{n+1}=\delta_{f} (_{0}h_{n}) \; u \; _{0}h_{n}              \; r \; _{0}h_{n}             \; d \; \delta_{f} (_{0}h_{n}) \\
_{3}h_{n+1}=\delta_{m}(_{0}h_{n})  \; u \; \delta_{g}(_{0}h_{n})  \; r \; \delta_{x}(_{0}h_{n}) \; d \; \delta_{m}(_{0}h_{n}) \\
_{4}h_{n+1}=\delta_{o} (_{0}h_{n}) \; u \; _{0}h_{n}              \; r \; _{0}h_{n}             \; d \; \delta_{f} (_{0}h_{n}) \\
_{5}h_{n+1}=\delta_{m}(_{0}h_{n})  \; u \; \delta_{g}(_{0}h_{n})  \; r \; \delta_{x}(_{0}h_{n}) \; d \; \delta_{x}(_{0}h_{n}) 
\end{array}
\end{equation}

and for the improper Hilbert curves

\begin{equation}
\begin{array}{l}
_{6}h_{n+1}=\delta_{f} (_{5}h_{n}) \; u \; \overline{ \delta_{m}(_{5}h_{n})} \; r \; _{5}h_{n} \; d \; \overline{ \delta_{y}(_{5}h_{n})} \\
_{7}h_{n+1}=\delta_{f} (_{5}h_{n}) \; u \; \overline{ \delta_{m}(_{5}h_{n})} \; r \; _{5}h_{n} \; d \; \delta_{a} (_{5}h_{n}) \\
_{8}h_{n+1}=\overline{\delta_{g} (_{5}h_{n})} \; u \; \overline{ \delta_{m}(_{5}h_{n})} \; r \; _{5}h_{n} \; d \; \delta_{a}(_{5}h_{n}) \\
_{9}h_{n+1}=\overline{\delta_{o} (_{5}h_{n})} \; u \;  \delta_{g}(_{5}h_{n}) \; r \; \overline{\delta_{a}(_{5}h_{n})} \; d \; \delta_{x} (_{5}h_{n}) \\
_{10}h_{n+1}=\delta_{m} (_{5}h_{n}) \; u \;  \delta_{g}(_{5}h_{n}) \; r \; \overline{\delta_{a}(_{5}h_{n})} \; d \; \overline{_{5}h_{n}} \\
_{11}h_{n+1}=\delta_{m} (_{5}h_{n}) \; u \;  \delta_{g}(_{5}h_{n}) \; r \; \overline{\delta_{a}(_{5}h_{n})} \; d \; \delta_{x}(_{5}h_{n}) 
\end{array}
\end{equation}

where

\begin{equation}
 \begin{array}{llll}
  \delta_{g}(u)=l & \delta_{g}(r)=u & \delta_{g}(d)=r & \delta_{g}(l)=d \\
  \delta_{x}(u)=r & \delta_{x}(r)=d & \delta_{x}(d)=l & \delta_{x}(l)=u \\
  \delta_{f}(u)=d & \delta_{f}(r)=l & \delta_{f}(d)=u & \delta_{f}(l)=r \\
  \delta_{m}(u)=d & \delta_{m}(r)=r & \delta_{m}(d)=u & \delta_{m}(l)=l \\
  \delta_{y}(u)=u & \delta_{y}(r)=l & \delta_{y}(d)=d & \delta_{y}(l)=r 
 \end{array}
\end{equation}

and the bar over the symbol represents the reversion operation.

\section{Conclusions}

Hilbert curves have proven to be useful constructions for a number of applications. It is often assumed, that there is exclusively one Hilbert curve in 2D, but this is not the case. There has been reports for five alternative patterns to which, in this contributions, six more patterns have been added by introducing the reversion operation. The twelve  patterns discussed in this paper have been called homogenous, to point out, that only one set of rules are applied to the $n$-order Hilbert curve, in order to generate the   $(n+1)$-order curve. All possible homogenous Hilbert curves in 2D have been exhausted. Is possible to construct additional Hilbert curves, by combining different set of rules at each quadrant (inhomogeneous curves) which will be subject of a future article.

The set of Hilbert curves here described, gives a larger flexibility in choosing the adequate scan pattern for each application. For example, close curves are useful in a path searching heuristic for the traveling salesman problem based on space filling curves. Other applications may benefit from curves that start and ends at opposite edges. Without regard of their possible use, the patterns which are here reported, can have different symmetries and structural features, which make each them interesting in their own right.

\section{Acknowledgments}

Universidad de la Habana and Universidad de la Ciencias Inform\'aticas are acknowledge for financial support and computational infrastructure.

We would like to thank an anonymous referee who made valuable suggestions that improved the final version of the manuscript. 


\begin{thebibliography}{10}
\expandafter\ifx\csname url\endcsname\relax
  \def\url#1{\texttt{#1}}\fi
\expandafter\ifx\csname urlprefix\endcsname\relax\def\urlprefix{URL }\fi
\expandafter\ifx\csname href\endcsname\relax
  \def\href#1#2{#2} \def\path#1{#1}\fi

\bibitem{sagan94}
H.~Sagan, Space-filling curves, Springer Verlag, New York, 1994.

\bibitem{chen05}
H.-L. Chen, Y.-I. Chang, Neighbor-finding based on space-filling curves, Inf.
  Syst. 30 (2005) 205--226.

\bibitem{chen11}
H.-L. Chen, Y.-I. Chang, All-nearest-neighbors finding based on the Hilbert
  curve, Expert Syst. with Appl. 38 (2011) 7462--7475.

\bibitem{bially69}
T.~Bially, Space-filling curves: Their generation and their application to
  bandwidth reduction, IEEE Trans. Inf. Th. IT-15 (1969) 658--664.

\bibitem{cox94}
J.~J. Cox, Y.~Takezaki, H.~R.~P. Ferguson, K.~E. Kohkonen, E.~L. Mulkay,
  Space-filling curves in tool-path applications., Computer-Aided Design 26
  (1994) 215--224.

\bibitem{songa02}
Z.~Songa, N.~Roussopoulosb, Using Hilbert curve in image storing and
  retrieving, Inform. Syst. 27 (2002) 523--536.

\bibitem{liang08}
J.-Y. Liang, C.-S. Chen, C.-H. Huang, L.~Liu, Lossless compression of medical
  images using Hilbert space-ﬁlling curves, Comp. Med. Img. and Graph. 32
  (2008) 174--182.

\bibitem{anders09}
S.~Anders, Visualization of genomic data with the Hilbert curve, Bioinformatics
  25 (2009) 1231--1235.

\bibitem{gao94}
J.~Gao, J.~M. Steele, General space filling curve heuristics and limit theory
  for the traveling salesman problem., J. Complex. 10 (1994) 230--245.

\bibitem{moon01}
B.~Moon, H.~V. Jagadish, C.~Faloutsos, Analysis of the clustering properties of
  the Hilbert space-filling curve, IEEE Trans. Knowl. Data. Eng. 13 (2001)
  124--141.

\bibitem{bauman06}
K.~E. Bauman, The dilation factor of the Peano-Hilbert curve, Math. Notes 80
  (2006) 609--620.

\bibitem{haverkort11}
H.~J. Haverkort, An inventory of three-dimensional Hilbert space-filling
  curves, CoRR abs/1109.2323.

\bibitem{moore00}
E.~H. Moore, On certain crinkly curves, Trans. Am. Math. Soc. 1 (1900) 72--90.

\bibitem{liu04}
L.~X, Four alternative patterns of the Hilbert curve, Appl. Math. Comp. 147
  (2004) 741--752.

\bibitem{bartholdi82}
J.~J. Bartholdi, L.~K. Platzman, An o(n log n ) planar traveling salesman
  heuristic based on space filling curves., Oper. Res. Lett. 1 (1982) 121--125.

\bibitem{seebold07}
P.~S\'e\'ebold, Tag system for the Hilbert curve, Discrete Math. and Th. Comp.
  Sci. 9 (2007) 213--226.

\end{thebibliography}

\section{References}

\pagebreak
\begin{figure}
\caption{The numeration of quadrants in clockwise orientation. Connectivity is given by $0 \rightarrow 1 \rightarrow 2 \rightarrow 3$.
}\label{fig:quadrants}
\end{figure}

\begin{figure}
\caption{Hilbert curve of order $n=1$. The solid circle marks the entry point, while the arrowhead the exit point of the curve.}\label{fig:h1}
\end{figure}

\begin{figure}
\caption{Hilbert curve of order $n=2$. a1: $_{0}H_{2}$, the original Hilbert curve; b1: $_{1}H_{2}$, the Moore curve; a2: $_{2}H_{2}$, Liu 1 curve; b2: $_{3}H_{2}$, Liu 2 curve; b2: $_{4}H_{2}$, Liu 3 curve; b2: $_{5}H_{2}$, Liu 4 curve.}\label{fig:h2}
\end{figure}

\begin{figure}
\caption{Affine transformation in each quadrant for the six proper Hilbert curves. The label at the lower right corner of each quadrant, follows the morphism of the tag system (see section \ref{sec:tag})}\label{fig:h2op}
\end{figure}

\begin{figure}
\caption{Proper Hilbert curve of order $n=3$. The labeling follows figure \ref{fig:h2}.}\label{fig:h3}
\end{figure}

\begin{figure}
\caption{Proper Hilbert curve of order $n=4$. The labeling follows figure \ref{fig:h2}.}\label{fig:h4}
\end{figure}

\begin{figure}
\caption{The reversion operation exchanges the entry and exit point of a curve.}\label{fig:swap}
\end{figure}

\begin{figure}
\caption{The affine transformations in each quadrant for the six improper Hilbert curves built from $_{5}H_{n}$. The line over the hands, flags that a reversion operation takes place. The label at the lower right corner of each quadrant, follows the morphism of the tag system (see section \ref{sec:tag}).}\label{fig:handimproper}
\end{figure}

\begin{figure}
\caption{The six improper Hilbert curves of order 3.}\label{fig:improper3}
\end{figure}

\begin{figure}
\caption{The six improper Hilbert curves of order 4.}\label{fig:improper4}
\end{figure}

\begin{figure}
\caption{The boundary vectors for the (a) $_{0}H_{2}$ curve and (b)  $_{5}H_{2}$ curve.}\label{fig:boundaryvector}
\end{figure}

\begin{figure}
\caption{The boundary vector diagrams for the homogeneous proper  Hilbert curves built from \ref{fig:boundaryvector}a. There are four possible orientation for the boundary vector in the first quadrant that complies with the connectivity to the second quadrant (labeled from 1 to 4).   Diagrams pairs 1b-3b and 2b-4b are each enantiomorphic pairs.}\label{fig:propercomplete}
\end{figure}

\begin{figure}
\caption{The boundary vector diagrams for the homogeneous improper  Hilbert curves built from \ref{fig:boundaryvector}b. There are four possible orientation for the boundary vector in the first quadrant that complies with the connectivity to the second quadrant (labeled from 1 to 4).  Diagrams pairs 1b-2b and 3b-4b are each enantiomorphic pairs. }\label{fig:impropercomplete}
\end{figure}

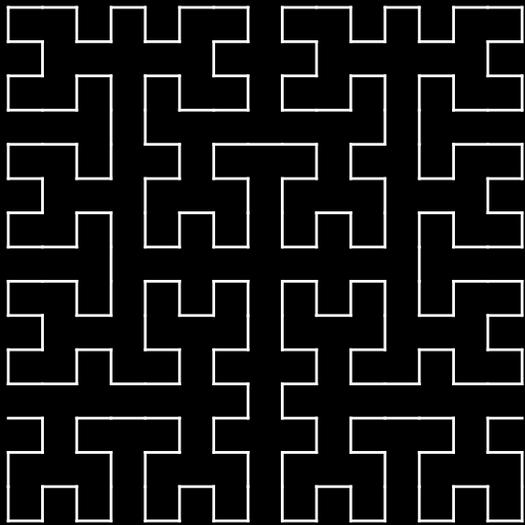
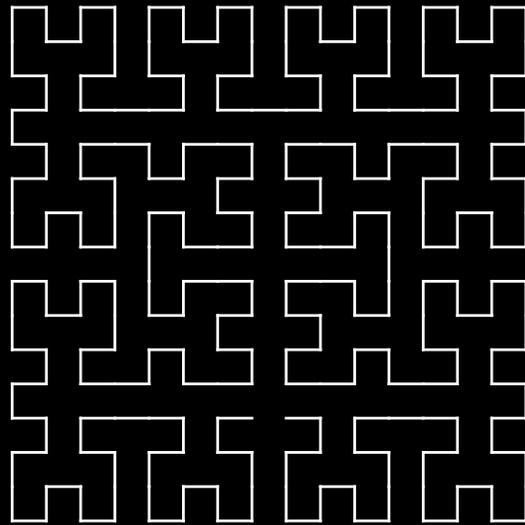
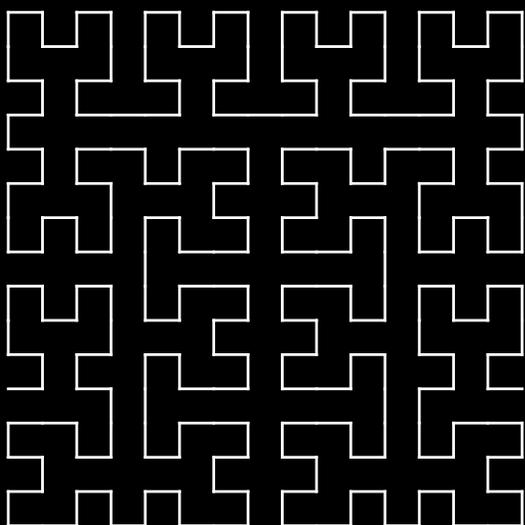
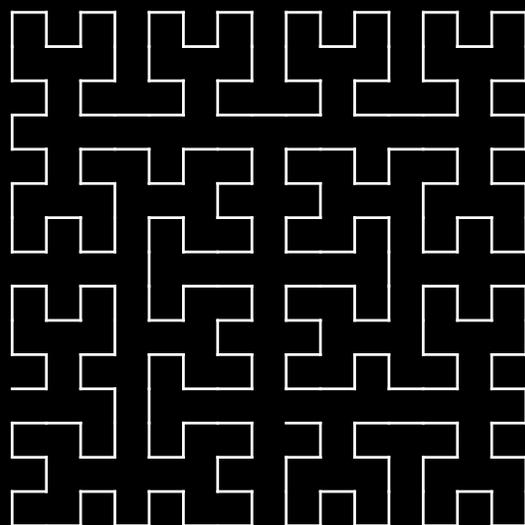
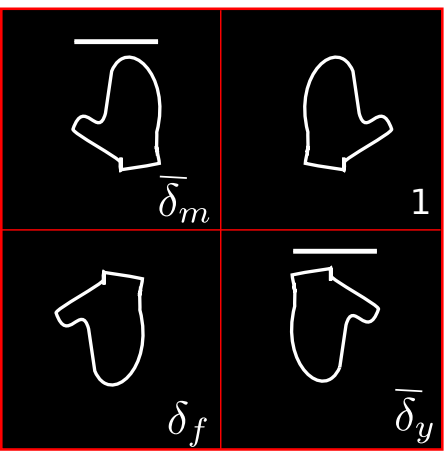
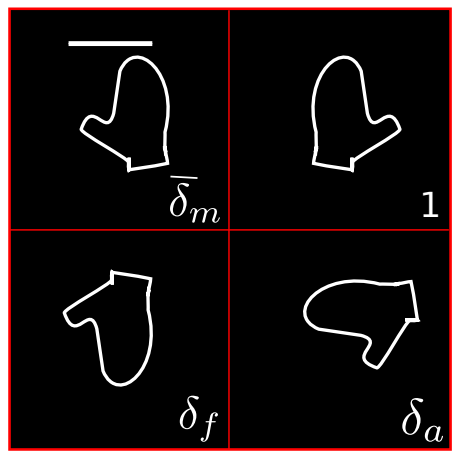
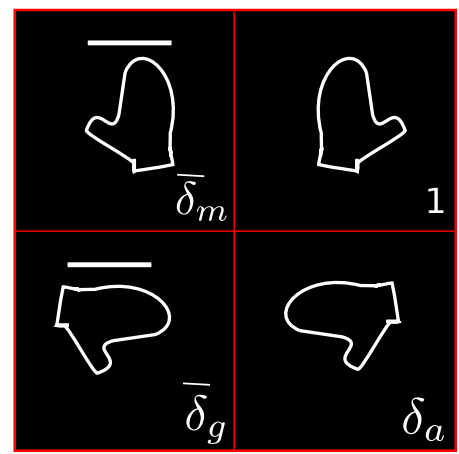
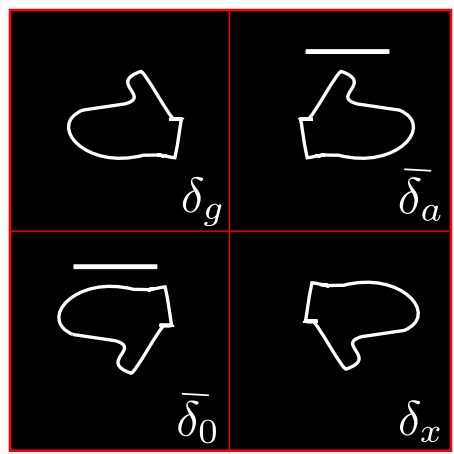
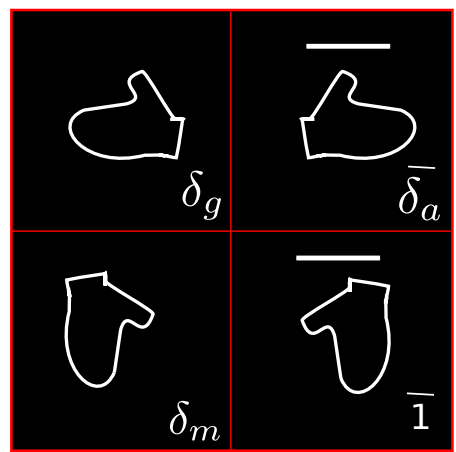
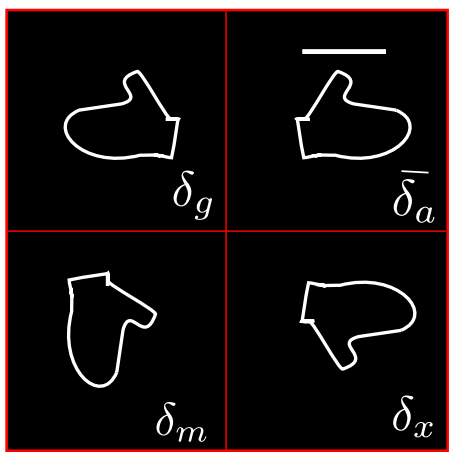
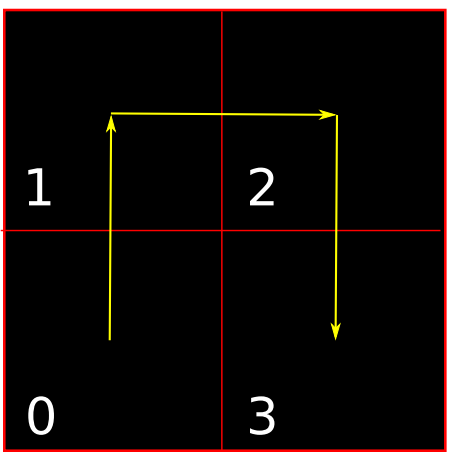
\end{document}